\documentclass[11pt]{article}

\usepackage{amsfonts}
\usepackage{amssymb}
\usepackage{amsbsy}
\usepackage{amsmath} 
 
\newcommand{\Matrix}[1]{\begin{bmatrix}
                 #1 
              \end{bmatrix}}

\newcommand{\bea}{\begin{array}}
\newcommand{\ea}{\end{array}}

\newtheorem{theorem}{Theorem}

\newtheorem{lemma}{Lemma}
\newtheorem{corollary}{Corollary}

\newtheorem{remark}{Remark}

\newenvironment{proof}{\begin{list}{$\!\!${\bf Proof.}%
  \rule{1pt}{0pt}}{\setlength{\leftmargin}{0pt}%
  \setlength{\itemindent}{30pt}%
  \setlength{\listparindent}{15pt}}\item}{\rule{0.3em}{0mm}%
  \hfill\framebox[1.2ex]{\rule{0.3em}{0mm}}\end{list}}

\newenvironment{prooft}[1]{\begin{list}{$\!\!${\bf Proof of Theorem~\ref{#1}.}%
  \rule{1pt}{0pt}}{\setlength{\leftmargin}{0pt}%
  \setlength{\itemindent}{30pt}%
  \setlength{\listparindent}{15pt}}\item}{\rule{0.3em}{0mm}%
  \hfill\framebox[1.2ex]{\rule{0.3em}{0mm}}\end{list}}

\newcommand{\Complex}{{\mathbb{C}}}
\newcommand{\re}{{\mathbb{R}}}

\newcommand{\Diag}{\mathop{\rm Diag}\nolimits}
\newcommand{\trace}{\mathop{\rm trace}\nolimits}
\newcommand{\real}{\mathop{\rm Re}\nolimits}
\newcommand{\RE}{\mathop{\rm Re}\nolimits}

\def\caa{{\mathcal A}}
\def\cbb{{\mathcal B}}
\def\ccc{{\mathcal C}}
\def\cdd{{\mathcal D}}
\def\cll{{\mathcal L}}
\def\cnn{{\mathcal N}}

\def\coo{{\mathcal O}}

\def\css{{\mathcal S}}
\def\mr#1{{(\ref{#1})}}

\newcommand{\beq}[2]{\begin{equation}
              \label{#1}
                {#2}
           \end{equation}}

\begin{document} 

\title{\textbf{A sequential semidefinite programming method and 
an application in passive reduced-order modeling}}

\author{Roland W.\ Freund\\
Department of Mathematics\\
   University of California at Davis\\
   One Shields Avenue\\
   Davis, California 95616, U.S.A.\\[2pt]
\textit{e-mail}: freund@math.ucdavis.edu\\
 \\ 
Florian Jarre and Christoph Vogelbusch\\
Institut f\"ur Mathematik \\
Universit\"at  D\"usseldorf\\
Universit\"atsstra{\ss}e 1 \\ 
D--40225 D\"usseldorf, Germany\\[2pt]
\textit{e-mail}: \{jarre,vogelbuc\}@opt.uni-duesseldorf.de}

\date{} 

\maketitle

\noindent
\textbf{Abstract.}
We consider the solution of nonlinear programs with nonlinear 
semidefiniteness constraints.
The need for an efficient exploitation of the cone of 
positive semidefinite matrices makes the solution of such
nonlinear semidefinite programs more complicated than the solution 
of standard nonlinear programs. 
In particular, a suitable symmetrization procedure needs to be
chosen for the linearization of the complementarity condition.
The choice of the symmetrization procedure
can be shifted in a very natural way to certain linear semidefinite 
subproblems, and can thus be reduced to a well-studied problem.
The resulting sequential semidefinite programming (SSP) method is a 
generalization of the well-known SQP method for standard nonlinear 
programs.
We present a sensitivity result for nonlinear
semidefinite programs,
and then based on this result, we give a self-contained 
proof of local quadratic convergence of the SSP method.
We also describe a class of nonlinear 
semidefinite programs that arise in passive reduced-order modeling,
and we report results of some numerical
experiments with the SSP method applied to problems in that class.

\bigskip

\noindent
\textbf{Key words.}
semidefinite programming, nonlinear programming,
sequential quadratic programming, 
quadratic semidefinite program,
sensitivity, convergence, reduced-order modeling, passivity, 
positive realness

\section{Introduction} \label{sec-intro}
In recent years, interior-point methods for solving linear 
semidefinite programs (SDPs) have received a lot of attention,
and as a result, these methods are now very well developed; see, e.g.,
\cite{todd,VanB96}, the papers in \cite{WolSV00}, and the references
given there.
At each iteration of an interior-point method, the complementarity 
condition is relaxed, symmetrized, and linearized.
Various symmetrization operators are known. 
The choice of the symmetrization operator and of the relaxation
parameter determine the step length at each iteration, and thus
the efficiency of the overall method.

In this paper, we are concerned with the solution of nonlinear semidefinite 
programs (NLSDPs).
Interior-point methods for linear SDPs can be extended to NLSDPs.
However, some additional difficulties arise.
First, the step length now also depends on the 
quality of the linearization of the nonlinear functions.
Second, the choice of the symmetrization procedure is considerably
more complicated than in the linear case since
the system matrix is no longer positive semidefinite.
To address these difficulties, in this paper, we consider an approach 
that separates the linearization and the symmetrization in a natural 
way, namely a generalization of the sequential quadratic
programming (SQP) method for standard nonlinear programs.
Such a generalization has already been mentioned by 
Robinson~\cite{robinson2} within the more general framework
of nonlinear programs over convex cones.
This framework includes NLSDPs as a special case. 
While Robinson did not discuss implementational issues of such 
a generalized SQP approach, the recent progress in the solution of
linear SDPs makes this approach especially
interesting for the solution of NLSDPs.
  
We first present a derivation of a generalized
SQP method, namely the sequential semidefinite programming (SSP) 
method, for solving NLSDPs.
In order to analyze the convergence of the SSP method,
we present a sensitivity result for certain local optimal solutions 
of general, possibly nonconvex, quadratic semidefinite programs. 
We then use this result to derive a self-contained proof of local quadratic 
convergence of the SSP method under the assumptions that the optimal 
solution is locally unique, strictly complementary, and satisfies a 
second-order sufficient condition.

One of the first numerical approaches for solving a class of NLSDPs 
was given in \cite{overton,overton92}.
Other recent approaches for solving NLSDPs are the
program package LOQO \cite{loqo} based on a 
primal-dual method; see also \cite{vanderbei}.
Another promising approach for solving large-scale SDPs
is the modified-barrier method proposed in \cite{kocvara}.
This modified-barrier approach does not require the  
barrier parameter to converge to zero, and may thus overcome some 
of the problems related to ill-conditioning in traditional 
interior-point methods.
Further approaches to solving NLSDPs have been presented
in \cite{ap1,ap2,ap3}. 
In \cite{ap1}, the augmented Lagrangian method is applied to NLSDPs,
while the approach proposed in \cite{ap2} 
is based on an SQP method generalized to NLSDPs.  
The paper \cite{ap2} also contains 
a proof of local quadratic convergence. 
However, in contrast to this paper, the algorithm \cite{ap2} 
is not derived from a comparison with interior-point algorithms,
and the proof of convergence does not use any differentiability
properties of the optimal solutions. 
In~\cite{ramirez}, Correa and Ram{\'{\i}}rez present
a proof of global convergence of a modification of the 
method proposed in \cite{ap2}. 
The modification employs certain merit functions to control 
the step lengths of the SQP algorithm. 
 
The remainder of this paper is organized as follows.
In Section~\ref{sec-notation}, we introduce some notation.
In Section~\ref{sec-rom}, we describe a class of nonlinear semidefinite
programs that arise in passive reduced-order modeling.
In Section~\ref{sec-linear}, we recall known results for 
linear SDPs in a form that can be 
easily transferred to NLSDPs.
In Section~\ref{sec-nonlinear}, we discuss primal-dual systems for
NLSDPs, and in Section~\ref{sec-SQP}, the SSP method is introduced
as a generalized SQP method.
In Section~\ref{sec-quad}, we present sensitivity results,
first for a certain class of quadratic SDPs,
and then for general NLSDPs.
Based on these sensitivity results, in Section~\ref{sec-conv}, we
give a self-contained proof of local quadratic convergence of the SSP method.
In Section\/~\ref{sec-num}, we present results of some numerical
experiments.
Finally, in Section~\ref{sec-conclusions}, we make some
concluding remarks.

\section{Notation} \label{sec-notation}

Throughout this article, all vectors  
and matrices are assumed to have real entries.
As usual, $Y^T = \Matrix{y_{ji}}$ denotes the
transpose of the matrix $Y = \Matrix{y_{ij}}$.
The vector norm $\|x\| := \sqrt{x^T x}$
is always the Euclidean norm and
$\|Y\| := \max_{\|x\|=1} \| Y x\|$ is the corresponding 
matrix norm.
For vectors $x\in \re^n$, $x\geq 0$ means that all entries of $x$
are nonnegative, and $\Diag(x)$ denotes the $n\times n$ diagonal
matrix the diagonal entries of which are the entries of $x$.
The $n\times n$ identity matrix is denoted by $I_n$.

The trace inner product on the space of real $n\times m$
matrices is given by
$$
\langle Z,Y \rangle := Z \bullet Y := \trace(Z^T Y)
 = \sum_{i=1}^n \sum_{j=1}^m z_{ij} y_{ij}
$$ 
for any pair $Y = \Matrix{y_{ij}}$ and $Z = \Matrix{z_{ij}}$
of $n\times m$ matrices.
The space of real symmetric $m\times m$ matrices is denoted 
by $\css^m$.
The notation $Y\succeq 0$ ($Y\succ 0$) 
is used to indicate that $Y\in\css^m$ is symmetric 
positive semidefinite (positive definite).

Semidefiniteness constraints are expressed by means of matrix-valued
functions from $\re^n$ to~$\css^m$.
We use the symbol $\caa: \re^n\to \css^m$ if such a function is
\textit{linear}, and the symbol $\cbb: \re^n \to \css^m$ 
if such a function is \textit{nonlinear}. 

Note that any linear function $\caa: \re^n \to \css^m$
can be expressed in the form
\beq{AA1}{
\caa(x) = \sum_{i=1}^n x_i A^{(i)}\quad \mbox{for all}\quad
x\in \re^n,
}
with symmetric matrices $A^{(i)}\in \css^m$, $i=1,2,\dots,n$.
Based on the representation~\mr{AA1} we introduce the norm
\beq{normA}{
\| \caa \| 
:= \biggl(\sum_{i=1}^n \bigl\|A^{(i)}\bigr\|^2 \biggr)^{\frac{1}{2}}
}
of $\caa$.
The adjoint operator $\caa^*:\css^m\to \re^n$ with 
respect to the trace inner product is defined by 
$$
\langle \caa(x) ,Y\rangle = \langle x, \caa^*(Y) \rangle
= x^T \caa^*(Y)
\quad \mbox{for all}\quad x\in \re^n \quad \mbox{and}\quad Y\in\css^m.
$$
It turns out that 
\beq{AA2}{
\caa^*(Y)=\Matrix{
           A^{(1)}\bullet Y \cr
           \vdots \cr
           A^{(n)}\bullet Y
         } \quad \mbox{for all}\quad Y\in\css^m.
}

We always assume that nonlinear functions $\cbb: \re^n \to \css^m$
are at least $\ccc^2$-differentiable.
We denote by
$$
B^{(i)}(x) := \frac{\partial}{\partial x_i}\cbb(x) \quad
\mbox{and} \quad
B^{(i,j)}(x):= 
    \frac{\partial^2}{\partial x_i\partial x_j}\cbb(x),\quad
i,j=1,2,\dots,n,
$$
the first and second partial derivatives of $\cbb$, respectively.
For each $x \in \re^n$, the derivative $D_x \cbb$ at $x$ induces
a linear function $D_x \cbb(x):\re^n \to \css^m$, which is given by
$$
D_{x} \cbb(x)[\Delta x] := 
\sum_{i=1}^n (\Delta x)_i B^{(i)}(x) \in \css^m 
\quad \mbox{for all}\quad \Delta x \in \re^n.
$$
In particular,
\[
\cbb(x + \Delta x) \approx \cbb(x) + D_{x} \cbb(x)[\Delta x],\quad
\Delta x \in \re^n, 
\]
is the linearization of $\cbb$ at the point $x$.
For any linear function $\caa: \re^n \to \css^m$, we have
\beq{AA1d}{
D_{x} \caa(x)[\Delta x] = \caa(\Delta x)
\quad \mbox{for all}\quad x,\, \Delta x \in \re^n.
}
We always use the expression on the right-hand side
of~\mr{AA1d} to describe derivatives of linear functions.

We remark that for any fixed matrix $Y \in\css^m$, the map
$x\mapsto\cbb(x)\bullet Y$ is a scalar-valued 
function of $x\in \re^n$.
Its gradient at $x$ is given by
\beq{AA3}{
\nabla_x \left(\cbb(x) \bullet Y \right) 
 = \bigl(D_x \left(\cbb(x) \bullet Y \right)\bigr)^T
         = \Matrix{
             B^{(1)}(x)\bullet Y \cr
             \vdots \cr
             B^{(n)}(x)\bullet Y
             } \in \re^n
}
and its Hessian by 
$$
\nabla^2_{x}\left( \cbb(x) \bullet Y \right) =
     \Matrix{
       B^{(1,1)}(x)\bullet Y & \cdots & B^{(1,n)}(x)\bullet Y\cr
        \vdots & & \vdots \cr
       B^{(n,1)}(x)\bullet Y &  \cdots & B^{(n,n)}(x)\bullet Y
       } \in \css^n.
$$
In particular, for any linear function $\caa: \re^n\to \css^m$,
in view of~\mr{AA1}, \mr{AA2}, and \mr{AA3}, we have
\beq{AA4}{
\nabla_x \left(\caa(x) \bullet Y \right) = \caa^*(Y).
}

\section{An application in passive reduced-order modeling} \label{sec-rom}

We remark that applications of linear SDPs
include relaxations of combinatorial optimization problems and 
problems related to Lyapunov functions or the positive real lemma 
in control theory; we refer the reader 
to \cite{Ali95,and,boydetal,freundjarre,todd,VanB96} and the
references given there.
In this section, we describe an application in passive reduced-order modeling
that leads to a class of NLSDPs. 

Roughly speaking, a system is called passive if it does not generate energy.
For the special case of time-invariant linear dynamical systems, passivity
is equivalent to positive realness of the frequency-domain transfer function
associated with the system.
More precisely, consider \textit{transfer functions} of the form
\beq{ZAA}{
Z_n(s) = B_2^T \bigl(G + s C\bigr)^{-1} B_1,\quad s\in \Complex,
}
where $G,\, C \in \re^{n\times n}$ and $B_1,\; B_2 \in \re^{n\times m}$
are given data matrices.
The integer $n$ is the \textit{state-space dimension} of
the time-invariant linear dynamical system, and $m$ is the number of
inputs and outputs of the system.
In~\mr{ZAA}, the matrix pencil $G + s C$ is assumed to be \textit{regular}, i.e.,
the matrix $G + s C$ is singular for only finitely many values of $s\in \Complex$.
Note that $Z_n$ is an $m\times m$-matrix-valued rational function of the 
complex variable $s\in \Complex$.

In \textit{reduced-order modeling}, one is given a large-scale 
time-invariant linear dynamical  
system of state-space dimension $N$, and the problem is to construct
a ``good'' approximation of that system of state-space dimension $n\ll N$;
see, e.g., \cite{freund} and the references given there.
If the large-scale system is passive, then for certain applications, it is
crucial that the
reduced-order model of state-space dimension $n$ preserves the passivity
of the original system.
Unfortunately, some of the most efficient reduced-order modeling techniques
do not preserve passivity.
However, the reduced-order models are often ``almost'' passive, and passivity
of the models can be enforced by perturbing the data matrices of the models.
Next, we describe how the problem of constructing such perturbations leads to
a class of NLSDPs.

An $m\times m$-matrix-valued rational function $Z$ is called \textit{positive real}
if the following three conditions are satisfied\/~$\! :$
\begin{itemize}
\item[\textup{(i)}]
$Z$ is analytic
in\/~$\Complex_+ := \{\, s\in \Complex \mid \RE(s) > 0\, \}$;
\item[\textup{(ii)}]
$Z(\overline{s}) = \overline{Z(s)}$ for all $s\in \Complex$;
\item[\textup{(iii)}]
$Z(s) + \bigl(Z(\overline{s})\bigr)^T \succeq 0$ for all $s\in \Complex_+$.
\end{itemize}

For functions $Z_n$ of the form~\mr{ZAA} positive realness (and thus
passivity of the system associated with $Z_n$) can be characterized
via linear SDPs; see, e.g.,~\cite{and,freundjarre} and the references given there.
More precisely, if the linear SDP
\beq{ZAC}{
\bea{rl}
P^T G + G^T P &\!\!\!\! \succeq 0,\\[4pt]
P^T C = C^T P &\!\!\!\! \succeq 0,\\[4pt]
       P^T B_1 &\!\!\!\! = B_2,
\ea
}
has a solution $P\in \re^{n\times n}$, then the transfer function~\mr{ZAA}, $Z_n$,
is positive real.
Conversely, under certain additional assumptions (see \cite{freundjarre}), 
positive realness of $Z_n$ implies the solvability of the linear SDP~\mr{ZAC}.

Now assume that $Z_n$ in~\mr{ZAA} is the transfer function of a 
non-passive reduced-oder model of a passive large-scale system.
Our goal is to perturb some of data matrices in~\mr{ZAA} so that the perturbed
transfer function is positive real.
For the special case $m=1$, such an approach is discussed in~\cite{bai}.
In this case, there is a simple eigenvalue-based characterization~\cite{bai0} 
of positive realness.
However, this characterization cannot be extended to the general case $m\geq 1$.
Another special case, which leads to linear SDPs, is described in \cite{phillips}.

In the general case $m\geq 1$, we employ perturbations $X_G$ and $X_C$ of
the matrices $G$ and $C$ in~\mr{ZAA}.
The resulting perturbed transfer function is then of the form
\beq{ZAF}{
\tilde{Z}_n(s) = B_2^T \bigl(G +X_G + s (C + X_C) \bigr)^{-1} B_1,
}
and the problem is to construct the perturbations $X_G$ and $X_C$  
such that $\tilde{Z}_n$ is positive real.
Applying the characterization~\mr{ZAC} of positive realness to~\mr{ZAF},
we obtain the following nonconvex NLDSP:
\beq{ZAG}{
\bea{rl}
P^T (G +X_G) + (G +X_G)^T P &\!\!\!\! \succeq 0,\\[4pt]
P^T (C + X_C) = (C + X_C)^T P &\!\!\!\! \succeq 0,\\[4pt]
       P^T B_1 &\!\!\!\! = B_2.
\ea
} 
Here, the unknowns are the matrices $P,\, X_G,\, X_C \in \re^{n\times n}$.
If~\mr{ZAG} has a solution $P,\, X_G,\, X_C$, then choosing the 
matrices $X_G$ and $X_C$ as the perturbations in~\mr{ZAF} guarantees
passivity of the reduced-order model given by the transfer 
function $\tilde{Z}_n$.

\section{Linear semidefinite programs} \label{sec-linear}

In this section, we briefly review the case of linear
semidefinite programs.

Given a linear function $\caa: \re^n \to \css^m$, a vector 
$b\in\re^n$, and a matrix $C\in\css^m$, a pair 
of primal and dual linear semidefinite programs is as follows:
\beq{plsdp}{
\bea{rl} 
\mbox{maximize} \quad C \bullet Y \quad \mbox{subject to} \quad
&\!\!\! Y \in \css^m,\ Y \succeq 0,\\[4pt]
&\!\!\!  \caa^*(Y) + b = 0,
\ea
}
and
\beq{dlsdp}{
\bea{rl}
\mbox{minimize} \quad b^T x \quad \mbox{subject to} \quad
&\!\!\! x \in \re^n,\\[4pt]
&\!\!\! \caa(x) + C \preceq 0.
\ea
}
We remark that this formulation is a slight variation of the
standard pair of primal-dual programs.
We chose the above version in order to facilitate the 
generalization of problems of the form \mr{dlsdp} to nonlinear 
semidefinite programs in standard form. 

If there exists a matrix $Y\succ 0$  that is feasible for \mr{plsdp},
then we  call $Y$ strictly feasible for \mr{plsdp} and say 
that \mr{plsdp} satisfies Slater's condition.
Likewise, if there exists a vector $x$ such that 
$\caa(x) + C \prec 0$, then we call \mr{dlsdp} strictly feasible and
say that \mr{dlsdp} satisfies Slater's condition.

The following optimality conditions for linear semidefinite
programs are well known;
see, e.g., \cite{shapiroetal}. 
If problem \mr{plsdp} or \mr{dlsdp} satisfies Slater's condition,
then the optimal values of \mr{plsdp} and \mr{dlsdp} coincide.
Furthermore, if in addition both problems are feasible,
then optimal solutions $Y^{\rm opt}$ and $x^{\rm opt}$ of 
both problems exist and $Y := Y^{\rm opt}$ and $x := x^{\rm opt}$ 
satisfy the complementarity condition 
\beq{opt1}{
Y S =0,\quad \mbox{where}\quad
S := -C- \caa(x). 
} 
Conversely, if $Y$ and $x$ are feasible points for \mr{plsdp} 
and \mr{dlsdp}, respectively, and satisfy the complementarity 
condition~\mr{opt1}, then $Y^{\rm opt} := Y$ is an optimal
solution of \mr{plsdp} and $x^{\rm opt} := x$ is an optimal solution 
of \mr{dlsdp}.

These optimality conditions can be summarized as follows.
If problem \mr{dlsdp} satisfies Slater's condition,
then for a point $x\in\re^n$ to be an optimal solution of \mr{dlsdp} 
it is necessary and sufficient that there exist
matrices $Y \succeq 0$ and $S\succeq 0$ such that 
\beq{sdp1}{
\bea{rl}
\caa(x) + C + S &\!\!\!\! = 0,\\[4pt]
\caa^*(Y) + b &\!\!\!\! = 0,\\[4pt]
          Y S &\!\!\!\! = 0.
\ea
}
Note that, in view of~\mr{AA4}, the second equation in~\mr{sdp1}
can also be written in the form
\beq{sdp2}{
\nabla_x\bigl((\caa(x)+C)\bullet Y\bigr) + b
 = \caa^*(Y) + b = 0.
}
Furthermore, the last equation in~\mr{sdp1} is equivalent to its
symmetric form, $Y S + S Y = 0$; see, e.g.,~\cite{aho}.
In the case of strict complementarity, the derivatives of $Y S=0$
and $Y S + S Y = 0$ are also equivalent.
For later use, we state these facts in the following lemma.

\begin{lemma}\label{lemma-facts}
Let\/~$Y,\, S \in \css^m$.
\begin{itemize}
\item[\textup{a)}] If\/~$Y\succeq 0$ or ~$S\succeq 0$, then
\beq{EQQ1}{
 Y S = 0\quad \Longleftrightarrow\quad  Y S + S Y = 0.
}
\item[\textup{b)}] If\/~$Y$ and $S$ are strictly complementary, 
 i.e.,\/~$Y,\, S \succeq 0$, $Y S = 0$, and\/~$Y + S \succ 0$, 
then  for any $\dot Y,\, \dot S \in \css^m$,
\beq{EQQ2}{
 Y \dot S + \dot Y S = 0\quad \Longleftrightarrow\quad  
 Y \dot S + \dot Y S + \dot S Y + S \dot Y = 0. 
}
Moreover,\/~$Y,\, S$ have representations of the form
\beq{EQQ3}{
Y = U \Matrix{Y_1 & 0\cr
                0 & 0} U^T, \quad 
S = U \Matrix{0 & 0\cr
              0 & S_2} U^T,
} 
where\/~$U$ is an $m\times m$ orthogonal matrix, $Y_1\succ 0$
is a $k\times k$ diagonal matrix, and  $S_2\succ 0$
is an $(m-k)\times (m-k)$ diagonal matrix, and any 
matrices $\dot Y,\, \dot S \in \css^m$ satisfying\/~$\mr{EQQ2}$ are of the form
\beq{EQQ4}{
\dot Y = U \Matrix{{\dot Y}_1 & {\dot Y}_3 \cr
                 {\dot Y}_3^T & 0} U^T,\quad               
\dot S = U \Matrix{0 & {\dot S}_3 \cr
        {\dot S}_3^T & {\dot S}_2} U^T,\quad \mbox{where}\quad
      Y_1 {\dot S}_3 + {\dot Y}_3 S_2 = 0. 
} 
\end{itemize} 
\end{lemma}

\begin{proof}
The equivalence~\mr{EQQ1} is well known; see, e.g.,~\cite[Page 749]{aho}.

We now turn to the proof of part b).
The strict complementarity of $Y$ and $S$ readily implies that $Y$ and $S$ 
have representations of the form~\mr{EQQ3}; 
see, e.g., \cite[Page 62]{luosturmzhang}.
Any matrices $\dot Y,\, \dot S \in \css^m$ can be written in the form
\beq{EQQ5}{
\dot Y = U \Matrix{{\dot Y}_1 & {\dot Y}_3 \cr
                 {\dot Y}_3^T & {\dot Y}_2} U^T,\quad               
\dot S = U \Matrix{{\dot S}_1 & {\dot S}_3 \cr
        {\dot S}_3^T & {\dot S}_2} U^T,
} 
where $U$ is the matrix from~\mr{EQQ3} and the block sizes in~\mr{EQQ5} are
the same as in~\mr{EQQ3}.
Using~\mr{EQQ3} and~\mr{EQQ5}, it follows that the equation on the left-hand
side of~\mr{EQQ2} is satisfied if, and only if,
\[
Y_1 {\dot S}_1 = 0,\quad 
{\dot Y}_2 S_2 = 0,\quad 
Y_1 {\dot S}_3 + {\dot Y}_3 S_2 = 0. 
\]
Since $Y_1$ and $S_2$ are in particular nonsingular, the first two relations
imply ${\dot S}_1 = 0$ and ${\dot Y}_2 =0$.
Thus, any matrices $\dot Y,\, \dot S \in \css^m$ satisfying the equation on
the left-hand side of \mr{EQQ2} are of the form~\mr{EQQ4}.
Similarly, using~\mr{EQQ3} and~\mr{EQQ5}, it follows that the equation on the 
right-hand
side of~\mr{EQQ2} is satisfied if, and only if,
\[ 
Y_1 {\dot S}_1 + {\dot S}_1 Y_1 = 0,\quad 
{\dot Y}_2 S_2 + S_2 {\dot Y}_2 = 0,\quad 
Y_1 {\dot S}_3 + {\dot Y}_3 S_2 = 0. 
\]
Since $Y_1 \succ 0$ and $S_2 \succ 0$, the first two relations
imply ${\dot S}_1 = 0$ and ${\dot Y}_2 =0$, and so $\dot Y,\, \dot S$ are again
of the form~\mr{EQQ4}.
\end{proof}
 
\section{Nonlinear semidefinite programs} \label{sec-nonlinear}

In this section, we consider nonlinear semidefinite programs,
which are extensions of the dual linear semidefinite 
programs~\mr{dlsdp}.

Given a vector $b\in \re^n$ and a matrix-valued function 
$\cbb:\re^n\to\css^m$, we consider
problems of the following form:
\beq{nlsdp}{
\bea{rl}
\mbox{minimize} \quad b^T x \quad \mbox{subject to} \quad
&\!\!\! x \in \re^n,\\[4pt]
&\!\!\!  \cbb(x) \preceq 0.
\ea
}
Here, the function $\cbb$ is nonlinear in general,
and thus~\mr{nlsdp} represents a class of nonlinear semidefinite 
programs. 
We assume that the function $\cbb$ is at least 
$\ccc^2$-differentiable.

For simplicity of presentation, we have
chosen a simple form of problem \mr{nlsdp}.
We stress that problem \mr{nlsdp} may also include additional 
nonlinear equality and inequality constraints. 
The corresponding modifications are detailed at the end
of this paper. 
Furthermore, the choice of the linear objective function $b^T x$
in~\mr{nlsdp} was made only to simplify notation.
A nonlinear objective function can always be 
transformed into a linear one by adding one artificial variable 
and one more constraint.
In particular, all statements about~\mr{nlsdp} in this paper
can be modified so that they apply to additional
nonlinear equality and inequality constraints and to 
nonlinear objective functions. 

Note that the class~\mr{nlsdp} reduces to linear semidefinite
programs of the form~\mr{dlsdp} if $\cbb$ is an affine
function.

The Lagrangian 
$\cll:\re^n \times \css^m  \to \re$ 
of \mr{nlsdp} is defined as follows:
\beq{Lag}{
\cll(x,Y) :=
b^T x + \cbb(x) \bullet Y.
}
Its gradient with respect to $x$ is given by
\beq{Lgrad}{
g(x,Y):=\nabla_x\cll(x,Y) =
b + \nabla_x\left(\cbb(x)\bullet Y \right) 
}
and its Hessian by
\beq{LHess}{
H(x,Y):=\nabla^2_x\cll(x,Y) =
 \nabla^2_x\left(\cbb(x)\bullet Y \right) .
}

If the problem \mr{nlsdp} is convex
and satisfies Slater's condition~\cite{Man69}, then for each
optimal solution $x$ of \mr{nlsdp} there
exists an $m\times m$ matrix $Y \succeq 0$
such that the pair $(x,Y)$ 
is a saddle point of the Lagrangian~\mr{Lag}, $\cll$. 

More generally, for nonconvex problems \mr{nlsdp}, let
$x\in\re^n$ be a feasible point of~\mr{nlsdp}, and
assume that the Robinson or Mangasarian-Fromovitz constraint 
qualification~\cite{Man69,robinson1,robinson2} is satisfied 
at $x$, i.e., there exists 
a vector $\Delta x\not=0$ such that $\cbb(x)+D_x \cbb(x)[\Delta x]
\prec 0$.
Then, if $x$ is a local minimizer of \mr{nlsdp}, the 
first-order optimality condition is satisfied, i.e., there exist
matrices $Y,\, S \in \css^m$ 
such that 
\beq{optcon}{
\bea{rl}
\cbb(x)+S &\!\!\!\! = 0,\\[4pt]
g(x,Y) &\!\!\!\! = 0,\\[4pt]
Y S & \!\!\!\! = 0,\\[4pt]
Y,\, S & \!\!\!\! \succeq 0.
\ea
}
The system~\mr{optcon} is a straightforward
generalization of the optimality conditions \mr{sdp1} and \mr{sdp2},
with the affine function $\caa(x) + C$ in~\mr{sdp2}
replaced by the nonlinear function $\cbb(x)$.

Primal-dual interior-point methods for solving \mr{nlsdp} roughly
proceed as follows.
For some sequence of duality parameters $\mu_k>0$,
$\mu_k\to0$, the solutions
of the perturbed primal-dual system,
\beq{pdipm}{
\bea{rl}
\cbb(x)+S& \!\!\!\! = 0,\\[4pt]
g(x,Y) &\!\!\!\! = 0,\\[4pt]
Y S & \!\!\!\! = \mu_k I_m,\\[4pt]
Y,\, S & \!\!\!\! \succ 0,
\ea
}
are approximated by some variant of Newton's method.
Since Newton's method does not preserve any inequalities,
the parameters $\mu_k>0$ are used to maintain strict
feasibility, i.e., $Y,\, S \succ 0$ 
for all iterates.

The solutions of \mr{pdipm} coincide with the solutions
of the standard logarithmic-barrier problems for \mr{nlsdp}.
Moreover, the logarithmic-barrier approach for solving \mr{nlsdp}
can be interpreted
as a certain choice of the `symmetrization operator'
for the equation, $Y S = \mu_k I_m$, in the  
third row of \mr{pdipm}; see Section \ref{sec-SQP} below.
With this choice, the barrier function yields a very natural 
criterion for the step-size control in trust-region algorithms.
The authors have implemented various versions of
predictor-corrector trust-region  
barrier methods for solving \mr{nlsdp}.
For a number of examples, the running times of the resulting
algorithms were comparable to the behavior
of interior-point methods for convex programs.
However, the authors also encountered several instances
in which the number of iterations for
these methods was very high compared to 
the typical number of iterations needed for solving linear SDPs.
For such negative examples it may be more efficient
to solve a sequence of linear SDPs in order to obtain an
approximate solution of \mr{nlsdp}. 
This observation motivated the SSP method described in the
next section.

\section{An SSP method for nonlinear semidefinite programs}
\label{sec-SQP}
 
In this section, we introduce the sequential semidefinite
programming (SSP) method, which is a generalization of the
SQP method for standard nonlinear programs to nonlinear
semidefinite programs of the form~\mr{nlsdp}.
For an overview of SQP methods for standard nonlinear programs, 
we refer the reader to~\cite{BogT95} and the references given 
there.

In analogy to the SQP method, at each iteration of the SSP method 
one solves a subproblem that is slightly more difficult than the 
linearization of \mr{optcon} at the current iterate.
More precisely, let $(x^k,Y^k)$ denote the current point
at the beginning of the $k$-th iteration.
One then determines corrections 
$(\Delta x, \Delta Y)$
and a matrix $S$  such that
\beq{sspproblem}{
\bea{rl}
 \cbb(x^k)+ D_x\cbb(x^k)[\Delta x]+S & \!\!\!\! = 0,\\[4pt]
\!\!\!\!\!\!
b + H^k \Delta x
+ \nabla_x \bigl(\cbb(x^k) \bullet (Y^k + \Delta Y)\bigr) 
 & \!\!\!\! = 0,\\[4pt]
(Y^k+\Delta Y^k) S & \!\!\!\! = 0,\\[4pt]
Y^k+\Delta Y,\, S & \!\!\!\! \succeq 0.
\ea
} 
Here and in the sequel, we use the notation
\beq{Hk}{
H^k := H\bigl(x^k,Y^k\bigr).
}
Recall from~\mr{Lgrad} and~\mr{LHess} that
$g(x,Y)$ and $H(x,Y)$ denote the gradient and Hessian
with respect to $x$, respectively, of the 
Lagrangian~\mr{Lag}, $\cll(x,Y)$, of the nonlinear
semidefinite program~\mr{nlsdp}.
Moreover, from~\mr{Lgrad} it follows that the
linearization of $g(x,Y)$ at the point $(x^k,Y^k)$
is given by 
\[
\bea{rl}
g\bigl(x^k + \Delta x,Y^k + \Delta Y\bigr)
 &\!\!\!\!  \approx b + H^k \Delta x 
 + \nabla_x \bigl(\cbb(x^k) \bullet (Y^k + \Delta Y)\bigr) .
\ea
\]
Thus the second equation in \mr{sspproblem} is just the
linearization of the second equation in \mr{optcon}.
Furthermore, the first equation of \mr{sspproblem}
is a straightforward linearization of the first equation in
\mr{optcon}.
This linearization is used in the same way in primal-dual 
interior methods.

The last two rows in~\mr{sspproblem} and \mr{optcon} are identical
when $Y$ in \mr{optcon} is rewritten as $Y = Y^k + \Delta Y$.
In analogy to SQP methods for standard nonlinear programs, the 
problem of how to guarantee the nonnegativity constraints,
namely $\cbb(x) \preceq 0$,
is thus shifted to the subproblem \mr{sspproblem}.
If the iterates $x^k$ generated from \mr{sspproblem} converge, then
their limit $x$ automatically satisfies $\cbb(x) \preceq 0$.

In contrast, interior methods use perturbations, symmetrizations, 
and linearizations for the last two rows in \mr{optcon},
resulting in cheaper linear subproblems that are typically 
less `powerful' than the subproblems \mr{sspproblem}.
An important aspect for interior methods for linear SDPs
is the choice of the symmetrization procedure for the 
bilinear equation $Y S = \mu_k I_m$; see, e.g., \cite{monteirozhang}.
For convex SDPs, 
theoretical convergence analyses are well developed
and also supported by numerical evidence of rapid 
convergence.
However, the generalization of these convergence results
to nonlinear nonconvex subproblems is far from obvious.
The proposed SSP method allows to apply the symmetrization to 
linear SDPs, thus reducing this aspect to a well-studied topic.

In summary, both the problem of choosing a suitable
symmetrization scheme and the problem of how to guarantee the 
nonnegativity constraints are shifted to the 
subproblem \mr{sspproblem}.

Note that the conditions \mr{sspproblem}
are the optimality conditions for the problem
\beq{sdpapprox}{
\bea{rl}
\mbox{minimize} \quad b^T \Delta x 
+ \frac{1}{2} (\Delta x)^T H^k \Delta x\quad \mbox{subject to}
\quad &\!\!\! \Delta x \in \re^n,\\[4pt]
  &\!\!\!  \cbb(x^k)+D_x\cbb(x^k)[\Delta x] \preceq 0.
\ea
}
The conditions \mr{sspproblem} and \mr{sdpapprox}
have been considered in \cite[Equations (2.1) and (2.2)]{robinson2}, 
with the remark that they have ``been found to be an appropriate
approximation of'' \mr{nlsdp} ``for numerical purposes''.
  
In order to be able to solve
the subproblem \mr{sspproblem} efficiently, in practice, one
replaces the matrix $H^k$ in \mr{sspproblem}, 
respectively~\mr{sdpapprox}, by a positive 
semidefinite approximation $\hat H^k$ of $H^k$. 
As in the case of standard SQP methods,
a BFGS update for the Hessian of the Lagrangian~\mr{Lag}, $\cll$,
can be used to approximate $H^k$
by some positive semidefinite matrix $\hat H^k$.
Given an estimate $\hat H^k$ of $H^k$ for the current, $k$-th,
SSP iteration, the quasi-Newton condition to generate a BFGS 
update $\hat H^{k+1}$
approximating the matrix $H^{k+1}$ for the next, $(k+1)$-st,
SSP iteration can be derived as follows:
\beq{approx.hess}{
\bea{rl}
\hat H^{k+1}\Delta x &\!\!\!\! = 
 \nabla_x \bigl(\cbb(x^{k+1})\bullet Y^{k+1}\bigr)
 - \nabla_x \bigl(\cbb(x^{k})\bullet Y^{k+1}\bigr)
    \\[4pt]
&\!\!\!\! =
   \nabla_x \cll\bigl(x^{k+1},Y^{k+1}\bigr)
 - \nabla_x \cll\bigl(x^k,Y^{k+1}\bigr)\\[4pt]
& \!\!\!\! \approx  \nabla^2_x \cll\bigl(x^{k+1},Y^{k+1}\bigr)
     \bigl(x^{k+1}-x^k\bigr).
\ea
}
If $\hat H^k$ is positive semidefinite, the BFGS  update
with the above condition can be suitably damped such that 
$\hat H^{k+1}$ is also positive semidefinite.
At each iteration of the SSP method, problem \mr{sdpapprox} is solved 
with $H^k$ replaced by the matrix $\hat H^{k+1}$
that is obtained by the BFGS update 
of $\hat H^k$ from the previous SSP iteration.
If $\hat H^{k+1}$ is positive semidefinite,
problem \mr{sdpapprox} essentially reduces to a linear SDP,
since the convex quadratic term in the objective function
can be written as a semidefiniteness constraint or 
a second-order-cone constraint. 
While the formulation as a second-order-cone constraint 
is more efficient,
and for example, can be specified as input for the program 
package SeDuMi \cite{sturm} in order to solve \mr{sdpapprox},
it was pointed out by \cite{noll} that
it may be most efficient to use a program that is designed for SDPs
with linear constraints and a convex quadratic objective function.

It seems that many results for standard SQP methods
carry over in a straightforward fashion to the SSP method.
For example, the SSP method 
can be augmented by a penalty term in case that 
the subproblems~\mr{sdpapprox} become infeasible. 
In this case,
the right-hand sides, ``$0$'', of the first three rows 
in \mr{sspproblem} are replaced by weaker, penalized 
right-hand sides.
Moreover, the convergence analysis of the method
proposed in \cite{ramirez,ap2} 
yields results that are comparable to the ones for
standard SQP methods. 

The standard analysis of quadratic convergence of SQP methods for 
nonlinear programs that satisfy strict complementarity
conditions proceeds by first showing that the 
active constraints will be identified correctly in the final
stages of the algorithm and then using the equivalence of the 
SQP iteration and the Newton iteration for the simplified KKT-system
in which only the active constraints are used.

For nonlinear semidefinite programs the situation is slightly more
complicated since 
it is difficult to identify active constraints.
The paper \cite{ap2} presents a proof that is based on a 
new approach by Bonnans et al. \cite{bonnans} and uses 
some general results due to Robinson \cite{robinson2}.
It does not require strict complementarity and
allows for approximate Hessian matrices in \mr{approx.hess}.

In the next two sections, we present a more elementary 
and self-contained approach to analyze convergence of the SSP 
method under a strict complementarity condition.

\section{Sensitivity results} \label{sec-quad}
 
In this section, we establish sensitivity results, first
for the special case of quadratic semidefinite programs
and then for general nonlinear semidefinite programs
of the form~\mr{nlsdp}.
More precisely, we show that strictly complementary solutions of 
such problems depend smoothly on the problem data.

We start with quadratic semidefinite programs of the form
\beq{nonl.sdp}{
\bea{rl}
\mbox{minimize} \quad  f(x) \quad \mbox{subject to}  \quad
&\!\!\! x \in \re^n,\\[4pt]
&\!\!\!  \caa(x) + C \preceq 0.
\ea
}
Here, $\caa:\re^n\to \css^m$ is a linear function,
$C \in \css^m$, and $f:\re^n \to \re$ is a quadratic function
defined by $f(x)=b^T x+ \frac{1}{2}x^T H x$, where
$b\in\re^n$ and $H\in\css^n$.
Note that we make no further assumptions on the matrix $H$.
Thus, problem~\mr{nonl.sdp} is a general, possibly nonconvex,
quadratic semidefinite program.

The problem \mr{nonl.sdp} is described by the
data
\beq{DATQ}{
\cdd := [\caa,b,C,H].
}
In Theorem~\ref{theorem1} below, we present a sensitivity
result for the solutions $\mr{nonl.sdp}$ when the data $\cdd$
is changed to $\cdd + \Delta \cdd$ where
\beq{DDATQ}{
\Delta \cdd := [\Delta \caa, \Delta b, \Delta C, \Delta H]
}
is a sufficiently small perturbation.
We use the norm
\[
\|\cdd\| :=  
\bigl( \|\caa\|^2 + \|b\|^2 + \|C\|^2 + \|H\|^2 \bigr)^{\frac{1}{2}}
\]
for data~\mr{DATQ} and perturbations~\mr{DDATQ}.
Recall that $\|\caa\|$ is defined in~\mr{normA}.

We denote by
$$
\cll^{\rm (q)}(x,Y) := f(x) + \left(\caa(x) + C\right) \bullet Y
$$ 
the Lagrangian of problem \mr{nonl.sdp}.
Note that $\nabla_x f(x) = b + H x$.
Together with~\mr{AA4}, it follows that
\[
\nabla_x \cll^{\rm (q)}(x,Y) = b + H x + \caa^*(Y)
\quad \mbox{and}\quad \nabla_x^2 \cll^{\rm (q)}(x,Y) = H.
\]

Recall that problem~\mr{nonl.sdp} is said to satisfy Slater's 
condition if there exists a vector $x \in \re^n$ 
with $\caa(x) + C \prec 0$.
Moreover, the triple $(\bar x,\bar Y,\bar S)$, where
$\bar x\in\re^n$ and $\bar Y,\, \bar S \in \css^m$, is
called a \textit{stationary point} of \mr{nonl.sdp} if 
\beq{strict.comp2}{
\bea{rl}
\caa(\bar x) + C + \bar S &\!\!\!\! = 0,\\[4pt]
b + H \bar x + \caa^*(\bar Y) &\!\!\!\! = 0,\\[4pt]
\bar Y \bar S  + \bar S \bar Y &\!\!\!\! = 0,\\[4pt]
\bar Y,\, \bar S & \!\!\!\! \succeq 0.
\ea
}
Here, we have used equivalence~\mr{EQQ1} of Lemma~\ref{lemma-facts}
and replaced $\bar Y \bar S = 0$ by its symmetric version,
which is stated as the third equation of~\mr{strict.comp2}.
If in addition to~\mr{strict.comp2}, one has
\beq{ZZA}{
\bar Y + \bar S\succ 0,
}
then $(\bar x,\bar Y,\bar S)$ is said to be a
\textit{strictly complementary} stationary point
of $\mr{nonl.sdp}$.

Let $\bar x\in\re^n$ be a feasible point of \mr{nonl.sdp}.
We say that $h\in\re^n$, $h \not= 0$, is 
a \textit{feasible direction} 
at $\bar x$ if $x= \bar x + \epsilon h$ is a feasible point 
of \mr{nonl.sdp} for all sufficiently small $\epsilon>0$.
Following \cite[Definition 2.1]{robinson2}, we say that
the \textit{second-order sufficient condition} holds at $\bar x$
with modulus $\mu>0$ if for all feasible directions $h\in \re^n$ 
at $\bar x$ with $h^T (b + H \bar x) = h^T \nabla_x f(\bar x) = 0$
one has
\beq{sos}{
h^T H h = h^T \bigl(\nabla_x^2 \cll^{\rm (q)}(\bar x, \bar Y)\bigr) h
 \ge \mu \|h\|^2.
}
After these preliminaries, our main result of this section 
can be stated as follows.

\begin{theorem} \label{theorem1}
Assume that problem $\mr{nonl.sdp}$
satisfies Slater's condition.
Let $(\bar x, \bar Y, \bar S)$ be a
locally unique and strictly complementary stationary point 
of $\mr{nonl.sdp}$ with data\/~$\mr{DATQ}$, $\cdd$, and assume that the
second-order sufficient condition holds at $\bar x$ with 
modulus $\mu>0$.
Then, for all sufficiently small perturbations\/~$\mr{DDATQ}$, 
$\Delta \cdd$, there exists a locally unique stationary point
$\bigl(\bar x(\Delta \cdd), \bar Y(\Delta \cdd), 
\bar S(\Delta \cdd)\bigr)$ of the perturbed program $\mr{nonl.sdp}$
with data\/~$\cdd + \Delta \cdd$.
Moreover, the point $\bigl(\bar x(\Delta \cdd), \bar Y(\Delta \cdd), 
\bar S(\Delta \cdd)\bigr)$ is a differentiable 
function of the perturbation $\mr{DDATQ}$, and for $\Delta \cdd = 0$,
$\bigl(\bar x(0), \bar Y(0), 
\bar S(0)\bigr) = (\bar x, \bar Y, \bar S)$.
The derivative $D_{\cdd}\bigl(\bar x(0),\bar Y(0),\bar S(0) \bigr)$
at $(\bar x,\bar Y,\bar S)$ is characterized by the directional
derivatives
$$
(\dot x,\dot Y,\dot S) := D_{\cdd}
\bigl(\bar x(0),\bar Y(0),
      \bar S(0) \bigr)[\Delta {\cdd}]
$$
for any $\Delta {\cdd}$.
Here $(\dot x,\dot Y,\dot S)$ is the unique solution
of the system of linear equations,
\beq{diff.eq2}{ 
\bea{rl}
\caa(\dot x) + \dot S &\!\!\!\!
=-\Delta C - \Delta \caa(\bar x), \\[4pt]
H \dot x + \caa^*(\dot Y) &\!\!\!\! = 
- \Delta b -\Delta H \bar x - \Delta \caa^*(\bar Y),\\[4pt] 
 \bar Y \dot S + \dot Y \bar S
 + \dot S \bar Y + \bar S \dot Y &\!\!\!\! =0,
\ea
}
for the unknowns $\dot x \in \re^n$, $\dot Y,\, \dot S \in \css^m$.
Finally, the second-order sufficient condition holds at
$\bar x(\Delta {\cdd})$ whenever $\Delta \cdd$ is sufficiently small.
\end{theorem}

\begin{remark} \label{rem1}
{\rm 
Theorem~\ref{theorem1} is an extension of
the sensitivity result for linear semidefinite programs
presented in~\cite{freundjarre2}.
A related sensitivity result for linear semidefinite programs
for a more restricted class of perturbations,
but also under weaker assumptions, is given in~\cite{sturmzhang}.
A local Lipschitz continuity property of unique and strictly
complementary solutions of linear semidefinite programs
is derived in~\cite{nayakk}.
}
\end{remark}

\begin{remark} \label{rem2}
{\rm 
While we did not explicitly state a linear independence
constraint qualification, commonly referred to as LICQ, it is
implied by our condition of uniqueness of the stationary
point; see, e.g.,~\cite{freundjarre2}.
Moreover, our assumptions on the stationary 
point $(\bar x, \bar Y, \bar S)$ imply that $\bar x$
is a strict local minimizer of \mr{nonl.sdp}.
}
\end{remark}

\begin{remark} \label{rem3}
{\rm 
The first and third equations in~\mr{diff.eq2} are symmetric
$m\times m$ matrix equations, and so only their upper triangular parts
have to be considered.
Thus the total number of scalar equations in~\mr{diff.eq2} is $m^2 + m + n$.
On the other hand, there are $m^2 + m + n$ unknowns, namely the entries
of $\dot x \in \real^n$ and of the upper triangular parts 
of $\dot Y,\, \dot S \in \css^m$.
Hence, \mr{diff.eq2} is a square system.
}
\end{remark}

\begin{remark} \label{rem4}
{\rm
In view of part b) of Lemma \ref{lemma-facts}, the 
last equation of \mr{diff.eq2} is equivalent to
\beq{diff.eq2b}{
 \bar Y\ \dot S + \dot Y\ \bar S = 0.
}
Thus, Theorem~\ref{theorem1} can be stated equivalently with~\mr{diff.eq2b}
in~\mr{diff.eq2}.
However, the resulting system of equations \mr{diff.eq2}
would then be overdetermined.
}
\end{remark}

\begin{prooft}{theorem1}
The proof is divided into four steps.

\noindent \textbf{Step 1.}
In this step, we establish the following result.
If the perturbed program has a local
solution that is a differentiable function of the perturbation,
then the derivative is indeed a solution of \mr{diff.eq2}.

Slater's condition 
is invariant under small perturbations of the problem data.
Hence, if there exists a local solution
$\bar x+\Delta x$, $\bar S+\Delta S$ 
of the perturbed problem near $\bar x$, $\bar S$, then
the necessary first-order conditions of the 
perturbed problem apply at $\bar x+\Delta x$, $\bar S+\Delta S$, 
and state that there exists a matrix $\Delta Y$ such that
$\bar Y+\Delta Y\succeq 0$, $\bar S+\Delta S\succeq 0$, and
\beq{pert3}{
\bea{rl}
(\caa+\Delta \caa)(\bar x + \Delta x)+ C + \Delta C + \bar S + \Delta S
          &\!\!\!\! = 0, \\[4pt]
    b + \Delta b + ( H + \Delta H) (\bar x + \Delta x)
      + (\caa^*+\Delta \caa^*) (\bar Y + \Delta Y) 
          &\!\!\!\! = 0, \\[4pt]
(\bar Y + \Delta Y)(\bar S + \Delta S)
 + (\bar S + \Delta S)(\bar Y + \Delta Y) &\!\!\!\! = 0.
\ea
}
Subtracting from these equations the first three equations of
\mr{strict.comp2} yields 
\beq{nleq}{
\bea{rl}
(\caa + \Delta \caa)(\Delta x) + \Delta S &\!\!\!\! =
 - \Delta C-\Delta \caa(\bar x), \\[4pt]
(H+\Delta H) \Delta x + (\caa^* + \Delta \caa^*)(\Delta Y) 
          &\!\!\!\! =
  - \Delta b -\Delta H \bar x - \Delta \caa^*(\bar Y),\\[4pt]
\bar Y\, \Delta S +  \Delta Y\, \bar S +  \Delta S\, \bar Y
  + \bar S\, \Delta Y  &\!\!\!\! 
   = - \Delta Y\, \Delta S - \Delta S\, \Delta Y.
\ea
}
Neglecting the second-order terms in \mr{nleq},
and using \mr{diff.eq2b},
we obtain the result claimed in \mr{diff.eq2}.
It still remains to verify the existence and differentiability
of $\Delta x$, $\Delta Y$, $\Delta S$. 

\noindent \textbf{Step 2.}
In this step, we prove that the system of linear equations \mr{diff.eq2}
has a unique solution.
To this end, we show that the homogeneous version of \mr{diff.eq2},
i.e., the system
\beq{scaled.syst2}{
\bea{rl}
{\caa}(\dot x) + \dot S &\!\!\!\! = 0,\\[4pt]
  H\dot x + {\caa}^*(\dot Y) 
     &\!\!\!\! = 0, \\[4pt]
 \bar Y \dot S + \dot Y \bar S
 + \dot S \bar Y + \bar S \dot Y &\!\!\!\! =0,
\ea
} 
only has the trivial solution $\dot x = 0$, $\dot Y = \dot S = 0$.

Let $\dot x \in \re^n$, $\dot Y,\, \dot S \in \css^m$ be any solution 
of~\mr{scaled.syst2}.
Recall that, in view of part b) by Lemma \ref{lemma-facts}
we may assume that $\bar Y$ and $\bar S$ are given in diagonal form:
\beq{diag}{
\bar Y = \Matrix{{\bar Y}_1 & 0\cr
                          0 & 0},\quad
\bar S  = \Matrix{0 & 0\cr
                  0 & {\bar S}_2},\quad
\mbox{where}\quad 
{\bar Y}_1,\, {\bar S}_2 \succ 0\quad \mbox{and}\quad
 {\bar Y}_1,\, {\bar S}_2\quad \mbox{are diagonal}.
}
Indeed, this can be done by replacing, in~\mr{scaled.syst2}, 
$\bar Y$, $\bar S$,
$\dot Y$, $\dot S$, $\caa(x)$ by $U^T \bar Y U$, $U^T \bar S U$,
$U^T \dot Y U $, $U^T \dot S U$, $U^T \caa(x) U$, respectively, 
where $U$ is the matrix in \mr{EQQ3},
and then multiplying the first and third rows from the left by $U$
and from the right by $U^T$.
Furthermore, in view of part b) by Lemma \ref{lemma-facts}, any
matrices $\dot Y\, \dot S \in \css^m$ satisfying the last equation 
of \mr{scaled.syst2} are then of the form
\beq{partition}{
\dot Y = \Matrix{{\dot Y}_1 & {\dot Y}_3 \cr
                 {\dot Y}_3^T & 0},\quad
\dot S = \Matrix{0 & {\dot S}_3 \cr
        {\dot S}_3^T & {\dot S}_2},\quad \mbox{where}\quad
      {\dot Y}_3 {\bar S}_2 + {\bar Y_1} {\dot S}_3 = 0.
}

Next, we establish the inequality 
\beq{dotx}{ 
{\dot x}^T H {\dot x} \geq \mu \| \dot x\|^2,
}
where $\mu>0$ is the modulus of the second-order sufficient
condition \mr{sos}.
Assume that $\dot x \not = 0$.
Let $\breve x \in \re^n$ be a Slater point for 
problem \mr{nonl.sdp}.
This guarantees that 
\beq{slat}{
M = \Matrix{M_{1} & M_{3} \cr
    M_{3}^T & M_{2}} := 
  -\bigl(\caa(\breve x) + C\bigr) \succ 0,
}
where the block partitioning $M$ is conforming with~\mr{partition}.
For $\eta > 0$, set
\beq{hath}{
h_{\eta}^+ : = \dot x + \eta (\breve x - \bar x)\quad \mbox{and}
\quad h_{\eta}^- : = -\dot x + \eta (\breve x - \bar x).
}
Since $\dot x \not = 0$, there exists an $\eta_0>0$ such
that $h_{\eta}^\pm \not= 0$ for all $0<\eta \leq \eta_0$.
Next, we prove that for all such $\eta$,
both vectors $h_{\eta}^+$ and $h_{\eta}^-$
are feasible directions for \mr{nonl.sdp} at $\bar x$.
Let $0<\eta \leq \eta_0$ be arbitrary, but fixed.
We then need to verify that 
$\caa(\bar x + \epsilon h_{\eta}^\pm) + C \preceq 0$ for 
all sufficiently small $\epsilon >0$.
Recall that $\caa$ is a linear function.
Using \mr{hath}, \mr{diag}, \mr{partition}, \mr{slat},
the first equation of \mr{strict.comp2}, and the first equation
of \mr{scaled.syst2},
one readily verifies that 
\beq{slat2}{
\bea{rl}
\caa(\bar x + \epsilon h_{\eta}^\pm) + C
&\!\!\!\! = (1 -\epsilon \eta) \bigl(\caa(\bar x) + C\bigr)
                 + \epsilon \eta  \bigl(\caa(\breve x) + C\bigr)
                 \pm \epsilon \caa(\dot x)\\[8pt]
&\!\!\!\! = -\left( \Matrix{ 0 & 0 \cr
                             0 & {\bar S}_2}
    + \epsilon \Matrix{ \eta M_{1} &  \eta M_{3} \pm {\dot S}_3 \cr
 (\eta M_{3} \pm {\dot S}_3)^T & 
     \eta (M_{2} - {\bar S}_2) \pm {\dot S}_2} \right).
\ea 
}
Recall that $\eta >0$ is fixed.
Since, by~\mr{diag} and~\mr{slat}, ${\bar S}_2 \succ 0$ 
and $M_1 \succ 0$, a standard Schur-complement argument shows
that the matrix on the right-hand side of \mr{slat2} is negative
definite for all sufficiently small $\epsilon >0$.
Thus the vectors~\mr{hath} 
are feasible directions for \mr{nonl.sdp} at $\bar x$
for any $\eta > 0$.
This in turn implies 
\beq{lss}{
{\dot x}^T (b+H\bar x) = {\dot x}^T \nabla_x f(\bar x) = 0.
}
Indeed, suppose that ${\dot x}^T \nabla_x f(\bar x) < 0$.
Then, for sufficiently small $\eta>0$, the feasible direction
$h_{\eta}^+$ also satisfies 
$\bigl(h_{\eta}^+\bigr)^T \nabla_x f(\bar x) < 0$,
and thus $h_{\eta}^+$ is a descent
direction for the objective function $f$ of \mr{nonl.sdp}
at the point $\bar x$.
This contradicts the local optimality of $\bar x$.
Likewise, if ${\dot x}^T \nabla_x f(\bar x) > 0$,
then, for sufficiently small $\eta>0$, 
$h_{\eta}^-$ is a descent
direction for the objective function $f$ of \mr{nonl.sdp}
at the point $\bar x$, leading to the same contradiction.
The second-order sufficient condition~\mr{sos} also holds
true on the closure of the feasible directions.
Since $\dot x$ is the limit of the feasible directions~\mr{hath}
for $\eta \to 0$ and $\dot x$ satisfies~\mr{lss}, the 
inequality~\mr{dotx} follows from~\mr{sos}.

Next recall from \mr{diag} that ${\bar Y}_1$ and ${\bar S}_2$ are 
positive definite diagonal matrices.  
The last relation in~\mr{partition} thus implies that
corresponding entries of the matrices ${\dot Y}_3$ and ${\dot S}_3$
are either zero or of opposite sign.
It follows that $\langle {\dot Y}_3, {\dot S}_3 \rangle \leq 0$, 
and equality holds if, and only if, ${\dot Y}_3 = {\dot S}_3 = 0$.
Using this inequality, together with the first two 
relations in~\mr{partition}, the first two equations of~\mr{scaled.syst2},
and~\mr{dotx}, one readily verifies that
\[ 
0 \geq 
2\, \langle {\dot Y}_3, {\dot S}_3 \rangle =
 \langle \dot Y, \dot S \rangle =
-\langle {\dot Y}, \caa(\dot x) \rangle = 
-\langle \caa^*(\dot Y) ,\dot x \rangle =
\langle H\dot x, \dot x \rangle = {\dot x}^T H {\dot x} \geq
\mu \| \dot x\|^2.
\]
Since $\mu > 0$, this implies
\beq{claim1}{ 
\dot x = 0 \quad \mbox{and} \quad {\dot Y}_3 = {\dot S}_3 = 0.
}
By the first row of~\mr{scaled.syst2}, it further follows that
\beq{claim2}{
{\dot S} = - \caa(\dot x) = - \caa(0) = 0.
}
Thus it only remains to show that $\dot Y = 0$.
In view of \mr{partition} and~\mr{claim2}, 
we have
\beq{claim3}{
\dot Y = \Matrix{{\dot Y}_1 & 0 \cr
                                    0 & 0}.
}
Now suppose that ${\dot Y}_1 \not = 0$.
Then, by~\mr{diag} and \mr{claim3}, we have 
\[
{\bar Y}_{\epsilon} := \bar Y + \epsilon \dot Y
\succeq 0\quad \mbox{and} \quad
{\bar Y}_{\epsilon}\not = \bar Y
\]
for all sufficiently small $|\epsilon|$.
Moreover, using~\mr{strict.comp2}, \mr{scaled.syst2}, and~\mr{claim2},
one readily verifies that the point
$(\bar x, {\bar Y}_{\epsilon}, \bar S)$ also 
satisfies~\mr{strict.comp2} for all sufficiently small $|\epsilon|$.
This contradicts the assumption that $(\bar x,\bar Y, \bar S)$
is a locally unique stationary point.
Hence ${\dot Y}_1 = 0$ and, by~\mr{claim3}, $\dot Y = 0$.

This concludes the proof that the square system \mr{scaled.syst2}
is nonsingular.

\noindent \textbf{Step 3.}
In this step, we show that the nonlinear system~\mr{strict.comp2} has a 
local solution that depends smoothly on the perturbation $\Delta \cdd$.
To this end, we apply the implicit function theorem to the system
\beq{red.eq}{
\caa(x) + C+ S = 0, \quad H x + b + \caa^*(Y) = 0,
\quad  Y\, S + S\, Y=0.
}
As we have just seen, the linearization 
of \mr{red.eq} at the point $(\bar x, \bar Y,  \bar S)$
is nonsingular, and hence \mr{red.eq}
has a differentiable and locally unique solution
$\bigl(\bar x(\Delta \cdd), \bar Y(\Delta \cdd),
\bar S(\Delta \cdd)\bigr)$.
Furthermore, we have $\bar Y(\Delta \cdd),\, \bar S(\Delta \cdd) \succeq 0$.
This semidefiniteness follows with a standard continuity argument:
The optimality conditions of the nonlinear SDP
coincide with the optimality conditions of the
linearized SDP. Under our assumptions, the latter
one has a unique optimal solution that depends
continuously on small perturbations of the
data; see, e.g., \cite{freundjarre2}.
Hence the linearized problem at the data
point $\cdd + \Delta \cdd$ has an optimal
solution $\bigl(\tilde x,\tilde Y, \tilde S\bigr)$
that satisfies the same optimality conditions as
$\bigl(\bar x(\Delta \cdd), \bar Y(\Delta \cdd),
\bar S(\Delta \cdd)\bigr)$. 
The solution of the linearized problem also satisfies $\tilde Y\succeq 0$,
$\tilde S\succeq 0$. 
Since $\bigl(\bar x(\Delta \cdd),
\bar Y(\Delta \cdd),\bar S(\Delta \cdd)\bigr)$ is
locally unique, it must coincide with
$\bigl(\tilde x,\tilde Y, \tilde S\bigr)$,
i.e., $\bar Y(\Delta \cdd),\, \bar S(\Delta \cdd)$ satisfy the semidefiniteness 
conditions.

\noindent \textbf{Step 4.}
In this step, we prove that the second-order sufficient
condition is satisfied at the perturbed solution.
Since feasible directions $h$ are defined only up to a positive
scalar factor, without loss of generality, one may require
that $\|h\| = 1$.
For the unperturbed problem~\mr{nonl.sdp},
the second-order sufficient condition at $\bar x$ then states
that $h^T H h \ge \mu $ for all $h\in\re^n$ with 
$\caa(x+\epsilon h) + C \preceq 0$, $\epsilon = \epsilon(h) >0$,
$\|h\|=1$, and $h^T(\bar b+H\bar x)=0$.
To prove that the second-order sufficient condition
is invariant under small perturbations $\Delta \cdd$ of the 
problem data $\cdd$, we thus need to show that for some fixed 
$\tilde \mu >0$, we have 
\beq{sosx1}{
h^T (H + \Delta H) h \ge \tilde \mu 
}
for all solutions $h\in \re^n$ of the
inequalities
\[
\bea{rl}
(\caa+\Delta \caa)\bigl(\bar x(\Delta \cdd) + \epsilon h\bigr) 
 + C + \Delta C \preceq 0, &\!\!\!\!
   \quad \epsilon = \epsilon(h) >0,\\[4pt]
 \| h\| = 1, \quad  &\!\!\!\!
h^T \bigl(b + \Delta b + (H + \Delta H)\bar x(\Delta \cdd)\bigr) = 0.
\ea
\]
In view of the first two relations in \mr{pert3}, the above is
equivalent to
\beq{FEAS}{
\epsilon (\caa+\Delta \caa)(h) \preceq \bar S(\Delta \cdd),\quad
\epsilon = \epsilon(h) >0,\quad  \| h\| = 1, \quad 
h^T (\caa^* + \Delta \caa^*) (\bar Y (\Delta \cdd)) = 0.
}
It remains to show that the set of solutions $h$ of~\mr{FEAS}
varies continuously with $\Delta D$.
Indeed, for any fixed $\tilde \mu$ with $0 < \tilde \mu < \mu$,
the second-order condition at $\bar x$ then readily 
implies that~\mr{sosx1} is satisfied for all solutions $h$ 
of~\mr{FEAS}, provided $\| \Delta \cdd\|$ is sufficiently small.

In Step 2, we have shown that both $\bar S(\Delta \cdd)$
and $\bar Y(\Delta \cdd)$
are continuous functions of $\Delta \cdd$ and that
the dimension of the null space of $\bar S(\Delta \cdd)$ 
is constant, namely equal to $k$, for all sufficiently 
small $\|\Delta D\|$. Moreover, the null space 
of $\bar S(\Delta \cdd)$ varies continuously with $\Delta D$.

Let $\Delta \cdd_k$ be a sequence of perturbations with
$\Delta \cdd_k\to 0$. 
Let $h_k$ be a sequence of associated solutions of 
\mr{FEAS}. 
It suffices to show that any accumulation point $\bar h$ of
the sequence $h_k$ satisfies \mr{FEAS} for $\Delta\cdd =0$
and the associated values $\Delta \caa =0$, 
$\bar Y(0)=\bar Y$, $\bar S(0)=\bar S$.
Since $\bar Y(\Delta\cdd)$ and $\Delta \caa^*$ vary continuously
with $\Delta \cdd$, it follows that $\bar h$ satisfies the
last two relations of \mr{FEAS} for $\Delta\cdd =0$.

We now assume by contradiction that 
$\epsilon\caa(\bar h)\not\preceq \bar S$
for any $\epsilon >0$. 
Since $\bar S\succeq 0$ this implies that there exists a 
vector $z \in \re^m$ 
with $\|z\|=1$, $z^T\caa(\bar h)z=\tilde \epsilon >0$,
and $z^T\bar Sz=0$. 
It follows that
\[
z^T(\caa +\Delta\caa_k)(h_k)z\ge \frac{\tilde \epsilon}{2}
\]  
if $k$ is sufficiently large.
Since the null space of $\bar S(\Delta\cdd)$
varies continuously with $\Delta D$, we have
\[
(z+\Delta z_k)^T\bar S(\Delta \cdd_k)(z+\Delta z_k)=0
\]
for some small $\Delta z_k\in \re^m$ whenever $\|\Delta \cdd_k\|$ is 
sufficiently small.
We now choose $\|\Delta \cdd_k\|$ so small, i.e., $k$ so large,
that 
\[
(z+\Delta z_k)^T(\caa +\Delta\caa_k)(h_k)(z+\Delta z_k)
\ge \frac{\tilde \epsilon}{4}.
\]
This implies that $h_k$ does not satisfy~\mr{FEAS}, and
thus yields the desired contradiction.
Hence $\bar h$ satisfies \mr{FEAS} for $\Delta\cdd =0$.
\end{prooft}

Theorem~\ref{theorem1} can be sharpened slightly.

\begin{corollary}\label{corollary1}
Under the assumptions of Theorem~\ref{theorem1} there exists a 
small neighborhood $\cnn$ of zero in the data space 
of $\mr{nonl.sdp}$ such that for all perturbations 
$\Delta \cdd \in \cnn$ of the problem data $\mr{DATQ}$, $\cdd$,
of $\mr{nonl.sdp}$,
there exists a local solution $(x_\Delta, Y_\Delta, S_\Delta)$ 
of $\mr{pert3}$
near $(\bar x, \bar Y, \bar S)$, at which the assumptions
of Theorem~\ref{theorem2} are also satisfied.
Moreover, the second derivatives
\[
\nabla^2_{\cdd} \bigl(x_\Delta,Y_\Delta,S_\Delta\bigr)[\Delta \cdd]
\]
of such local solutions $(x_\Delta, Y_\Delta, S_\Delta)$
are uniformly bounded for all $\Delta \cdd \in \cnn$.
\end{corollary}

\begin{proof}
The first part of the corollary is an immediate consequence of 
Theorem~\ref{theorem1}. 
For the second part observe that
the second derivative is obtained by differentiating the
system \mr{diff.eq2}.
For sufficiently small perturbations $\Delta \cdd$,
the singular values of this system are uniformly bounded
away from zero, and hence the second derivatives 
are uniformly bounded.
\end{proof}

Theorem~\ref{theorem1} can be generalized to the
class of NLSDPs of the form~\mr{nlsdp}.
Recall that, by~\mr{Lag} and~\mr{LHess}, the
Lagrangian of~\mr{nlsdp} and its Hessian are given by
\beq{Lag2}{
\cll(x,Y) = b^T x + \cbb(x) \bullet Y\quad \mbox{and}\quad
H(x,Y) = \nabla^2_x\left(\cbb(x)\bullet Y \right),
}
respectively.
The generalization of Theorem~\ref{theorem1} to 
problems~\mr{nlsdp} is then as follows.

\begin{theorem} \label{theorem2}
Let $x^*$ be a local
solution of~$\mr{nlsdp}$, and let $Y^*$ be an
associated Lagrange multiplier.
Assume that the Robinson constraint 
qualification is satisfied at $x^*$
and that the point $(x^*,Y^*)$
is strictly complementary and locally unique.
Finally, assume that the second-order sufficient condition holds
at $x^*$ with modulus $\mu>0$. 
Then, $\mr{nlsdp}$
has a locally unique solution for small perturbations
of the data $(\cbb,b)$, and the solution depends smoothly on the 
perturbations.
\end{theorem}

\begin{proof}
First, we define the linear function 
$\caa := D_x \cbb(x^*): \re^n \to \css^m$, and the
matrices $C:=\cbb(x^*)$ and $H:=H(x^*,Y^*)$.
Then, the SSP approximation \mr{sdpapprox} (with $p=q=0$)
of \mr{nlsdp} at the point $(x^*,Y^*)$ is 
just the quadratic semidefinite problem
\beq{simple.sqp.prob}{
\bea{rl}
\mbox{minimize} \quad  b^T\Delta x
               + \frac{1}{2} (\Delta x)^T H \Delta x
\quad \mbox{subject to}  \quad
&\!\!\! \Delta x \in \re^n,\\[4pt]
&\!\!\!  \caa(\Delta x) + C \preceq 0.
\ea
}
Note that~\mr{simple.sqp.prob} is a problem of the 
form~\mr{nonl.sdp} with data \mr{DATQ}, $\cdd$.
Moreover, $\Delta \bar x:=0$, $\bar Y:=Y^*$, 
and $\bar S:= - \caa(0) - C$ satisfy the 
conditions \mr{strict.comp2}. 
These conditions coincide with the first-order conditions 
of~\mr{nlsdp}, and thus the point 
$(\Delta \bar x,\bar Y, \bar S)$ is also the unique 
solution of \mr{strict.comp2}. 
Furthermore, the second-order sufficient condition 
for \mr{simple.sqp.prob} and~\mr{nlsdp} coincide. 
This condition guarantees that $\Delta\bar x$ is a locally unique
solution of~\mr{simple.sqp.prob}.
Finally, the Robinson constraint qualification 
for problem~\mr{nlsdp} at $x^*$
implies that problem \mr{simple.sqp.prob} satisfies Slater's condition.
In particular, all assumptions of Theorem~\ref{theorem1} are 
satisfied. 
Small perturbations $\Delta \cdd$ of the data 
of~\mr{nlsdp} result in small changes of the
corresponding SSP problem~\mr{simple.sqp.prob}.
Since Theorem~\ref{theorem1} allows for arbitrary changes in
all of the data of~\mr{simple.sqp.prob}, the claims follow.
\end{proof} 

\section{Convergence of the SSP method} \label{sec-conv}

In this section, we prove that the plain SSP method with 
step size 1 is locally quadratically convergent.

For pairs $(x,Y)$, where $x\in \re^n$, $Y\in \css^m$,
we use the norm
\[
\bigl\| (x,Y) \bigr\|
:= \bigl( \|x\|^2 + \|Y\|^2 \bigr)^{\frac{1}{2}}.
\]
The main result of this section can then be stated as
follows.

\begin{theorem}\label{theorem3}
Assume that the function $\cbb$ 
in $\mr{nlsdp}$ is $\ccc^3$-differentiable
and that problem $\mr{nlsdp}$ has a locally unique and
strictly complementary solution $(\bar x,\bar Y)$ 
that satisfies the Robinson constraint qualification and the
second-order sufficient condition with modulus $\mu>0$, cf. $\mr{sos}$.
Let some iterate $(x^k,Y^k)$ be given
and let the next iterate 
$$
(x^{k+1},Y^{k+1})
:=(x^k,Y^k)+ (\Delta x,\Delta Y)
$$ 
be defined as the local solution of $\mr{sspproblem}$, or, equivalently,
$\mr{sdpapprox}$, that is closest to $(x^k,Y^k)$.
Then there exist $\epsilon >0$ and $\gamma <1/\epsilon$ such that 
$$
\bigl\| (x^{k+1},Y^{k+1})-
   (\bar x,\bar Y) \bigr\| \, \le \,
\gamma \bigl\| (x^{k},Y^{k})-
 (\bar x,\bar Y) \bigr\|^2
$$ 
whenever $\bigl\| (x^{k},Y^{k})-
(\bar x,\bar Y)\bigr\| < \epsilon$.
\end{theorem}

\begin{proof}
The proof is divided into three steps.
In the first step, we establish the exact relation of
problems \mr{sdpapprox} and \mr{nonl.sdp}.
In a second step, we consider a point $x^k$ near $\bar x$. 
We show that $x^k$ is the optimal solution of an SSP subproblem
the data of which is at most $\coo(\|x^{k}-\bar x\|)$ 
away from the data of the SSP subproblem at $(\bar x, \bar Y)$.
We remark that $x^{k}$ is always the optimal solution of the 
$(k-1)$-st subproblem, but the data of this subproblem
lies $\coo(\|x^{k-1}-\bar x\|)$ away from the SSP subproblem 
at $(\bar x, \bar Y)$.
In a third step, we then show by a perturbation
analysis that the correction $\Delta x=x^{k+1}-x^k$ is of size
$\coo(\|x^{k}-\bar x\|+\|Y^{k}-\bar Y\|)$ and
that the residual for the SSP subproblem in 
the $(k+1)$-st step is of size $\coo((\|x^{k}-\bar x\|+\|Y^{k}-\bar Y\|)^2)$.

\noindent \textbf{Step 1.}
We first show how the SSP subproblem \mr{sdpapprox} 
can be written in the form \mr{nonl.sdp}.
To this end, we define the linear function
$ \caa:=D_x\cbb(x^k):\re^n\to\css^{m}$, and the
matrices $C:=\cbb(x^k)$ and $H := H(x^k,Y^k)$.
Note that the linear constraint $\caa(\Delta x)+C\preceq 0$ is just the
linearization of the nonlinear constraint 
$\cbb(x^k+\Delta x) \preceq 0$ about the point $x^k$.
Finally, let $b$ be as in \mr{nlsdp}. 
The SSP subproblem \mr{sdpapprox} 
then takes the simple form \mr{simple.sqp.prob}, and in particular,
it conforms with the format \mr{nonl.sdp} of Theorem~\ref{theorem2}.

\noindent \textbf{Step 2.}
Let any point $x^{k}$ close to $\bar x$ be given.
We show that $\Delta x =0$ is a 
local solution of a problem of the form \mr{simple.sqp.prob},
where the data is `close' to the data of \mr{sdpapprox}
at $\left(\bar x,\bar Y\right)$.
Let $\Delta C:=\cbb(\bar x)-\cbb(x^k)$.
By continuity of $\cbb$, $\|\Delta C\|$ is
of order $\coo(\|x^k-\bar x\|)$. 
Let
\[ 
\Delta b:= - \nabla_x\left(\cbb(x^k)\bullet \bar Y\right)- b
=  -\caa^*(\bar Y)  - b.
\] 
From \mr{Lgrad} and the second row of \mr{optcon}, we obtain
the estimate 
$\|\Delta b\|=\coo(\|x^k-\bar x\|)$,
and the point $(0,\bar Y,\bar S)$ satisfies the first-order 
conditions,
\[
\caa(0)+C+\Delta C+ \bar S=0, \quad 
b + \Delta b + H \cdot 0 + \caa^*(\bar Y) = 0,\quad 
\quad \bar Y\, \bar S=0, 
\]
for the quadratic semidefinite program
\beq{approxk}{
\bea{rl}
\mbox{minimize} \quad  (b+\Delta b)^T\Delta x + 
   \frac{1}{2} (\Delta x)^T H \Delta x \quad  
\mbox{subject to} \quad  &\!\!\! \Delta x \in \re^n,\\[4pt]
&\!\!\! \caa(\Delta x) + C+\Delta C \preceq 0.
\ea
}
Let 
\beq{opt.data}{ 
\bar \caa:= D_x\cbb(\bar x),\quad \bar b := b,
\quad \bar C := \cbb(\bar x),
\quad \bar H:= \nabla^2_x \cll(\bar x,\bar Y),
}
be the data of the SSP subproblem \mr{sdpapprox} at the 
point $(\bar x,\bar Y)$.
Then, the data of \mr{approxk} differs from the data
\mr{opt.data} by a perturbation of norm
$\coo(\|x^k-\bar x\|+\|Y^k-\bar Y\|)$.
Here, the term $\|Y^k-\bar Y\|$ reflects the fact that $H$ 
and $\bar H$ also differ by the choice of $Y$.
Note that the point $(\Delta x,Y,S) = (0,\bar Y,\bar S)$ 
satisfies the first-order optimality conditions,
\beq{opt.residual}{
\bar \caa(\Delta x)+ \bar C+  S = 0, \quad
b + \bar H \Delta x + \bar \caa^*( Y) = 0,\quad
Y S=0,
} 
of the quadratic problem \mr{simple.sqp.prob} with data \mr{opt.data}.
As shown in Theorem~\ref{theorem2}, the assumptions for
the nonlinear SDP \mr{nlsdp} at $(\bar x,\bar Y)$ imply
that problem \mr{simple.sqp.prob} with 
data \mr{opt.data} satisfies all conditions of Theorem~\ref{theorem2}
at $(0,\bar Y,\bar S)$.

Since the second-order sufficient condition depends
continuously on the data of \mr{nonl.sdp},
it follows that for \mr{approxk}, the second-order condition
at $(0,\bar Y,\bar S)$ is satisfied, provided that
$\|x^k-\bar x\|$ and $\|Y^k-\bar Y\|$
are sufficiently small.
Thus problem~\mr{approxk} fulfills all assumptions of
Theorem~\ref{theorem2}.

\noindent \textbf{Step 3.}
By definition, $(\Delta x,Y)=(0,\bar Y)$ 
is the optimal solution (with associated multiplier) of \mr{approxk}.
Let $(x^k,Y^k)$ be close to $(\bar x,\bar Y)$.
The SSP subproblem replaces the data $\Delta b$ and $\Delta C$
of \mr{approxk} by 0 (of the respective dimension).
Thus, the data of \mr{approxk} is changed by a perturbation 
of order $\coo(\|x^k-\bar x\|+\|Y^k-\bar Y\|)$.
We assume that this perturbation lies in the neighborhood $\cnn$
about zero as guaranteed by Corollary~\ref{corollary1}.
Denote the optimal solution of the SSP subproblem by
$(\Delta x,\bar Y+\Delta \hat Y)$.

The SSP subproblem is then used to define $(x^{k+1},Y^{k+1})$. Let 
\[
\caa^+:=D_x\cbb(x^{k+1}),\quad 
C^+:=\cbb(x^{k+1}),\quad
H^+:=\nabla^2_x \cll(x^{k+1},Y^{k+1})
\]
be the data of the SSP subproblem at the next, $(k+1)$-st iteration.

Corollary~\ref{corollary1} states that 
$(\Delta x,\Delta \hat Y)$  
are given by the tangent equations \mr{diff.eq2} plus some uniformly 
bounded second-order terms. 
Thus, $(\Delta x,\Delta \hat Y)$ are of the order 
$\coo(\|x^k-\bar x\|+\|Y^k-\bar Y\|)$.
Here, $\Delta \hat Y$ is a correction of the unknown point $\bar Y$,
while the correction $\Delta Y=Y^{k+1}-Y^k$ produced by the SSP subproblem has 
the form $\Delta Y = \Delta \hat Y + Y^k - \bar Y$. 
Obviously, also the norm $\|\Delta Y\|$ of this correction 
is of the order $\coo(\|x^k-\bar x\|+\|Y^k-\bar Y\|)$.

Next, we compute an upper bound on the size of the 
residuals of the first and second equations
in \mr{opt.residual} at $(x^{k+1},Y^{k+1},S^{k+1})$.
Note that the residual term of the third equation 
in \mr{opt.residual} is zero.
By definition of $(\Delta x, Y^{k+1},S^{k+1})$, it follows that
\beq{zero.residual}{
\caa(\Delta x) +C+ S^{k+1} = 0, \quad
b + H \Delta x + \caa^*(Y^{k+1}) = 0,\quad
Y^{k+1}\, S^{k+1}=0.
}
If the data of \mr{nlsdp} is $C^3$-smooth, this implies that
\[
\bea{rl}
(\caa^+)^*(Y^{k+1})+b
&\!\!\!\! = \caa^*(Y^{k+1}) +b+ \Delta \caa^*(Y^{k+1}) \\[4pt]
&\!\!\!\! = -H \Delta x+ \Delta \caa^*(Y^{k})
            + \Delta \caa^*(\Delta Y)\\[4pt]
&\!\!\!\! = - \nabla^2_x \bigl(\cbb(x^k)\bullet Y^{k}\bigr)\, \Delta x
   +\bigl(\nabla_x\cbb(x^k+\Delta x) 
   -\nabla_x\cbb(x^k)\bigr)\bullet Y^k 
            + \Delta \caa^*(\Delta Y)\\[4pt]
&\!\!\!\! = \coo(\|\Delta x\|^2+\|\Delta Y\|^2),
\ea
\]
where $\Delta\caa:=\caa^+-\caa$.
Likewise, it follows from \mr{zero.residual} that 
\[
\bea{rl}
C^++S^{k+1} &\!\!\!\! =\Delta C+C+S^{k+1}
                      =\Delta C-\caa(\Delta x) \\[4pt]
&\!\!\!\! = \cbb(x^{k+1})-\cbb(x^{k})-D_x\cbb(x^k)[\Delta x]
          = \coo(\|\Delta x\|^2),
\ea
\]
where $\Delta C:= C^+-C$. 
Hence, we can define perturbations $\hat b$ and $\hat C^+$
of $b$ and $C^+$ of order
$\|\hat b-b\|+\|\hat C^+-C^+\|= \coo((\|x^k-\bar x\|+\|Y^k-\bar Y\|)^2)$
such that 
$(\Delta x,Y,S)=(0,Y^{k+1},S^{k+1})$ is an optimal solution
of the problem~\mr{simple.sqp.prob} 
with data $\caa^+$, $\hat b$, $\hat C^+$, $H^+$.
By the same derivation as above, the next SSP step
has length $\coo((\|x^k-\bar x\|+\|Y^k-\bar Y\|)^2)$
and generates residuals of 
order $\coo((\|x^k-\bar x\|+\|Y^k-\bar Y\|)^4)$.
Repeating this process, it then follows by a standard argument that 
$\|x^{k+1}-\bar x\|$ and $\|Y^{k+1}-\bar Y\|$
are of order $\coo((\|x^k-\bar x\|+\|Y^k-\bar Y\|)^2)$ as well. 
\end{proof}

\begin{remark} \label{rem5}
{\rm 
As mentioned before, one will typically choose to solve
SSP subproblems with a positive semidefinite approximation $\hat H$
to the Hessian of the Lagrangian. 
A proof of convergence for such modifications
is the subject of current research; see, e.g.,~\cite{ramirez}.
Since all the data enters in a continuous fashion in the preceding
analysis, it follows that the SSP method with step size one 
is still locally superlinearly convergent if the matrices 
$H^k$ in \mr{Hk} are replaced by approximations $\hat H^k$
with $\|H^k-\hat H^k\|\to 0$.
}
\end{remark}

\begin{remark} \label{rem6}
{\rm 
The assumption of $\ccc^3$-differentiability of the function
$\cbb$ in Theorem~\ref{theorem3} can be weakened to
$\ccc^2$-differentiability and a Lipschitz condition for the Hessian
at $\bar x$.
}
\end{remark}

The result of Theorem~\ref{theorem3} can be extended
to the following slightly more general class of NLSDPs.
Given a vector $b\in \re^n$, a matrix-valued function 
$\cbb:\re^n\to\css^m$, and two vector-valued functions
$c:\re^n\to\re^p$ and $d:\re^n\to\re^q$, we consider
problems of the following form:
\beq{nlsdpg}{
\bea{rl}
\mbox{minimize} \quad b^T x \quad \mbox{subject to} \quad
&\!\!\! x \in \re^n,\\[4pt]
&\!\!\!  \cbb(x) \preceq 0, \\[4pt]
&\!\!\! c(x)  \le 0, \\[4pt]  
&\!\!\! d(x)  =  0.
\ea
}
The Lagrangian of problem \mr{nlsdp} takes the form
$\cll:\re^n \times \css^m \times \re^p \times \re^q \to \re$:
\beq{Lagg}{
\cll(x,Y,u,v) :=
b^T x + \cbb(x) \bullet Y + u^T c(x) + v^T d(x).
}
Its gradient with respect to $x$ is given by
\beq{Lgradg}{
g(x,Y,u,v):=\nabla_x\cll(x,Y,u,v) =
b + \nabla_x\left(\cbb(x)\bullet Y \right) +
   \nabla_x c(x)\, u +  \nabla_x d(x)\, v
}
and its Hessian by
\beq{LHessg}{
H(x,Y,u,v):=\nabla^2_x\cll(x,Y,u,v) =
 \nabla^2_x\left(\cbb(x)\bullet Y \right) +
\sum_{i=1}^p u_i \nabla_x^2 c_i(x)
+ \sum_{j=1}^q v_j \nabla_x^2 d_j(x).
}
Note that in~\mr{Lgradg}, the gradients of the vector-valued functions
$c(x)$ and $d(x)$ are defined as $\nabla_x c(x) := (D_x c(x))^T$
and $\nabla_x d(x) := (D_x d(x))^T$.

For NLSDPs~\mr{nlsdpg}, the SSP subproblems are of the form 
\beq{sdpapproxg}{
\bea{rl}
\mbox{minimize} \quad b^T \Delta x 
+ \frac{1}{2} (\Delta x)^T H^k \Delta x\quad \mbox{subject to}
\quad &\!\!\! \Delta x \in \re^n,\\[4pt]
  &\!\!\!  \cbb(x^k)+D_x\cbb(x^k)[\Delta x] \preceq 0, \\[4pt]
  &\!\!\!  c(x^k) + D_x c(x^k)\, \Delta x \le 0,\\[4pt]
  &\!\!\!  d(x^k) + D_x d(x^k)\, \Delta x = 0.
\ea
}
The extension of Theorem \ref{theorem3} to problems \mr{nlsdpg}
is as follows.

\begin{theorem}\label{theorem4}
Assume that the functions $\cbb$, $c$, and $d$ in \mr{nlsdpg}
are $\ccc^3$-differentiable,
and that problem $\mr{nlsdpg}$ has a locally unique and
strictly complementary solution $(\bar x,\bar Y, \bar u, \bar v)$
that satisfies the Robinson constraint qualification and the
second-order sufficient condition with modulus $\mu>0$, cf. $\mr{sos}$.
Let some iterate $(x^k,Y^k,u^k,v^k)$ be given
and let the next iterate
$$
(x^{k+1},Y^{k+1},u^{k+1},v^{k+1})
:=(x^k,Y^k,u^k,v^k)+ (\Delta x,\Delta Y,\Delta u,\Delta v)
$$
be defined as the local solution of $\mr{sdpapproxg}$ that is closest 
to $(x^k,Y^k,u^k,v^k)$.
Then there exist $\epsilon >0$ and $\gamma <1/\epsilon$ such that
$$
\bigl\| (x^{k+1},Y^{k+1},u^{k+1},v^{k+1})-
   (\bar x,\bar Y, \bar u, \bar v) \bigr\| \, \le \,
\gamma \bigl\| (x^{k},Y^{k},u^{k},v^{k})-
 (\bar x,\bar Y, \bar u, \bar v) \bigr\|^2
$$
whenever $\bigl\| (x^{k},Y^{k},u^{k},v^k)-
(\bar x,\bar Y, \bar u, \bar v)\bigr\| < \epsilon$.
\end{theorem}

\begin{proof}
By our assumption on strict complementarity, all entries of
the vector $\bar v$ 
of the Lagrange multipliers associated with the
equality constraints $d(x)=0$ of \mr{nlsdp} are different
from zero.
Without loss of generality, we assume that $\bar v > 0$.
Indeed, for any entry $\bar v_j <0$ we replace 
the corresponding constraint $d_j(x)=0$ by the equivalent constraint 
$-d_j(x)=0$. 
These sign changes do not change the iterates generated
by \mr{sdpapprox}.
Moreover, for $\bigl(x^{k},Y^{k},u^{k},v^{k}\bigr)$ 
sufficiently close to $\bigl(\bar x,\bar Y, \bar u, \bar v\bigr)$ 
it follows from $\bar v>0$ that the iterates do not change
when the constraints  $d(x)=0$ are replaced by $d(x)\le 0$.
We can thus assume that $q=0$, i.e., there are no equality 
constraints in \mr{nlsdpg}.

We further assume that, without loss of generality,
the matrix $\cbb$ is augmented to a $2\times 2$ block diagonal 
matrix, where the $(2,2)$-block is the diagonal matrix $\Diag(c(x))$.
Thus, for the analysis of the SSP method we may assume 
that $p=q=0$ in \mr{nlsdp}, i.e., we only need to consider
problems of the form \mr{nlsdp}.
\end{proof}
   
\section{Numerical results} \label{sec-num}

In this section, we present results of some numerical experiments with a
Matlab implementation of the SSP method.
Actually, our Matlab program is for a slightly more general class of 
nonlinear programs with conic constraints (NLCPs).
The numerical experiments with our program illustrate the theoretical results 
of the preceding sections.
In particular, quadratic convergence is observed for problems where
the Hessian $H$ of the Lagrangian at the optimal solution 
is positive semidefinite. 
In cases where $H$ is not positive semidefinite, our implementation uses
perturbations of the nonconvex SSP subproblems in order to obtain
convex conic subproblems. 
In these cases, typically, the rate of convergence of the algorithm
based on such perturbed problems is only linear.

The Matlab program generates its search directions by solving conic
quadratic subproblems using Version 1.05R5 of SeDuMi~\cite{sturm}. 
SeDuMi allows free and positive
variables as well as Lorentz-cone (``ice-cream cone'') constraints, rotated
Lorentz-cone constraints, and semidefinite cone constraints. 
The NLCPs can also be formulated in terms of these cones.
In order to simplify the use of SeDuMi for the SSP subproblems,
the NLCPs are rewritten in the following standard format: 
\beq{sedumiformatnlcp}{
\bea{rl}
\mbox{minimize} \quad c^T x \quad \mbox{subject to} \quad
&\!\!\! x\in K,\\[4pt]
&\!\!\! F(x) = 0.
\ea
}
Here, $K$ is a Cartesian product of free variables and several
cones of the types allowed in SeDuMi.

We tested the following techniques
for generating positive semidefinite approximations of $H$:
a BFGS approach, the Hessian of the augmented Lagrangian, and 
the orthogonal projection of $H$ onto the cone of positive semidefinite matrices. 
Our experiences with these techniques are as follows.
\begin{enumerate}
\item
The BFGS approach can result in considerably more SSP iterations compared 
to the projection of the Hessian of the Lagrangian. 
Moreover, the BFGS approach strongly depends on the
initial matrix $H_0$.
A good choice is $H_0:=V \max(D,\epsilon I) V^T$, where $H = V D V^T$ is the
eigenvalue decomposition of $H$.
\item 
The use of the Hessian of the augmented
Lagrangian can be a good choice for some problems, but for most of our test
problems the penalty parameter had to be very large to obtain a semidefinite
Hessian. 
This, in turn, significantly reduced the precision of our computations.
\item
In spite of not being affinely invariant, the use of the 
projection of the Hessian of the Lagrangian resulted in the 
most efficient overall algorithm.
\end{enumerate}

We also tested different step length strategies.
\begin{enumerate}
\item
The following penalty line search with a
quadratic correction gave good results for all test cases.
The SSP subproblem
provides a search direction $\Delta x$ for problem \mr{sedumiformatnlcp}.
By solving a least-squares problem, a vector $q$ is
computed satisfying $D_x F(x)q=-F(x+\Delta x)$.
For $\lambda\in[0,1]$, a line search along
the points $x(\lambda):= x+\lambda\Delta x+\lambda^2 q$
is performed based on the penalty function
$$
M \| F(x(\lambda))\| + c^Tx(\lambda),
$$
where $M>0$ is a penalty parameter.
\item
For some examples, the choice of a filter approach was slightly better. 
In the filter approach used here, a Euclidean 
trust-region radius was always set to be 
1.5 times larger than the previous step, and 
non-acceptable steps were not discarded, but instead an Armijo-type
step-length reduction was used to generate an acceptable step. 
The motivation for this modified filter strategy lies in
the fact that the computation of a solution of a subproblem
is very expensive, and 
therefore discarding the solution of a subproblem is avoided. 
The above filter approach led to very fast
convergence, especially for convex problems. 
\item
For the examples presented here, the trust-region approach was the best
choice. 
The SSP subproblem was restricted by an additional 
Euclidean trust-region constraint.
For problems of the form \mr{roland_bsp} below, it was sufficient to apply
the trust-region constraint only to the variables $X_G,\, X_C$, while $P,\, S$ 
remain free. 
For these examples, an
additional corrector step significantly accelerated the convergence.
For this corrector step, $X_G,\, X_C$ is kept fixed, and $P,\, S$ are updated
by solving an additional linear SDP. 
At each iteration, the ratio between predicted and actual 
reduction was computed.
Depending on that ratio, the  step was accepted and the trust region was 
increased or decreased, or the step was rejected and the trust 
region was decreased.
\end{enumerate}

For our numerical examples, we use nonlinear nonconvex SDPs of the form~\mr{ZAG},
which we rewrite in the form~\mr{roland_bsp} below.
Recall that in~\mr{ZAG}, $G,\, C \in \re^{n\times n}$ and 
$B_1,\; B_2 \in \re^{n\times m}$ are given data matrices.
The nonconvex NLSDP used for the numerical examples is then as follows:
\beq{roland_bsp}{
\bea{rl}
\mbox{minimize} \quad  \| S\| \quad \mbox{subject to} \quad
&\!\!\! P\in \re^{n\times n},\ S\in \re^{n\times m},\\[4pt]
&\!\!\! X_G\in \re^{n\times n},\ \Vert X_G\Vert \le r_G,\\[4pt]
&\!\!\! X_C\in \re^{n\times n},\ \Vert X_C\Vert \le r_C,\\[4pt]
&\!\!\! P^T B_1+S = B_2, \\[4pt]
&\!\!\! P^T(G+X_G)+(G+X_G)^T P  \succeq  \varepsilon_G I,\\[4pt]
&\!\!\! P^T(C+X_C)+(C+X_C)^T P  \succeq  \varepsilon_C I,\\[4pt]
&\!\!\! P^T(C+X_C)-(C+X_C)^T P  =  0.
\ea
}
Furthermore, in~\mr{roland_bsp}, in addition to the constraints on the norms of 
the perturbations $X_G$ and $X_C$, we restrict $X_G$ and $X_C$ to have possible
nonzero entries only in certain positions, which depend on the nonzero structure
of the given matrices $G$ and $C$, respectively.
For our numerical tests, the data matrices $G$ and $C$ in~\mr{roland_bsp}
are generated as follows.
First, two matrices $C_{org}$ and $G_{org}$ were constructed such that the 
associated transfer function is guaranteed to be positive real.
Then, certain entries of $G_{org}$ and $C_{org}$ were perturbed
by random perturbations of norm at most $\varepsilon_G$ and $\varepsilon_C$
respectively, to define the data matrices $G$ and $C$.
In all our examples, the transfer functions of the systems given by the resulting matrices
$C$ and $G$ were not positive real.

All our computations were run on a Xeon with a clock rate of 2.8 GHz 
and 3 GB RAM.
All solutions were computed to a precision of 12 decimal digits. 
In the following table, we list the problem dimension $n$,
the total number $M(n)$ of equality constraints,
the total number $N(n)$ of scalar unknowns, the number of iterations, and
the cpu time (in seconds).

\bigskip

\begin{center}
\begin{tabular}{|r|r|r|r|r|} \hline 
$n$ & $M(n)$ & $N(n)$ & iterations & cpu time\\
\hline \hline 
8&  118 & 285& 5& 3.71 \\
9&  146 & 348& 5& 4.43 \\
10& 177 & 417& 7& 8.06 \\
11& 211 & 492& 5& 7.18 \\
12& 248 & 573& 8& 16.05 \\
13& 288 & 660& 4& 10.88 \\
14& 331 & 753& 6& 20.77 \\
15& 377 & 852& 7& 30.12 \\
16& 426 & 957& 6& 34.38 \\
17& 478 & 1068& 5& 37.40 \\
18& 533 & 1185& 10& 91.17 \\
19& 591 & 1308& 4& 47.61 \\
20& 652 & 1437& 5& 83.66 \\
21& 716 & 1572& 4& 289.48 \\
22& 783 & 1713& 4& 416.31 \\
23& 853 & 1862& 5& 151.33 \\
24& 926 & 2013& 6& 683.54 \\
25& 1002 & 2176& 3& 145.60 \\
26& 1081 & 2337& 5& 612.22 \\
27& 1163 &  2508& 7& 518.92 \\
28& 1248&  2685& 5& 789.41 \\
29& 1336&  2868& 4& 475.52 \\
30& 1427&  3057& 7& 4213.50 \\
31& 1521&  3252& 4& 784.34 \\
32& 1618&  3455& 6& 4659.64 \\
33& 1718&  3660& 5& 1130.44 \\ 
34& 1821& 3877& 2& 630.53 \\
35& 1927& 4092& 6& 1799.36 \\
\hline 
\end{tabular}
\label{tab:pade}
\end{center}

\bigskip

Table \ref{tab:pade} shows that the number of iterations is nearly
independent of the dimension $n$ of the problem, while---as expected---the cpu time increases 
with~$n$. 
The total number of constraints is approximately $M(n) \approx \frac{3}{2}n^2$,
and the total number of scalar variables 
is approximately $N(n) \approx 3n^2$. 
The number of iterations to solve the linear semidefinite subproblems not only depends on the
dimension, but also on other properties of the
problem as, for example, a comparison of the problems of dimension 32 and 33 shows. 
In this case, the iteration counts differ only by one, yet the cpu time quadruples, since
SeDuMi needs more iterations to solve the subproblems. 
Some of the linear semidefinite subproblems are nearly infeasible, a situation for which 
SeDuMi (and other solvers) needs a higher number of interior-point iterations.
 
\section{Concluding remarks} \label{sec-conclusions}

We have discussed the SSP method, which is a generalization of 
the SQP method for standard nonlinear programs to nonlinear 
semidefinite programming problems.
For the derivation of this generalization, we
have chosen a motivation that contrasts the SSP method
with primal-dual interior methods. 
For interior methods that are applied directly to nonlinear 
semidefinite programs, the choice of the symmetrization procedure 
is considerably more complicated than in the linear case since
the system matrix is no longer positive semidefinite.
In the proposed method, the choice of the symmetrization
scheme is shifted in a very natural way to the subproblems, 
and is thus reduced to a well-studied problem.
Our convergence analysis differs from the convergence analyses
of standard SQP methods in that it is based on a sensitivity result 
for certain optimal solutions of quadratic semidefinite programs.
The derivation of this sensitivity result is also of
independent interest.

\end{document}